\documentclass[12pt,oneside,english]{amsart}
\usepackage[T1]{fontenc}
\usepackage[latin9]{inputenc}
\usepackage{amsthm}
\usepackage{amstext}
\usepackage{amssymb}

\makeatletter

\newcommand{\noun}[1]{\textsc{#1}}
\providecommand{\tabularnewline}{\\}

\numberwithin{equation}{section}
\numberwithin{figure}{section}
  \theoremstyle{definition}
  \newtheorem*{example*}{\protect\examplename}

\usepackage{bussproofs,hyperref}
\usepackage{tikz}
\usepackage{tikz-cd}
\usetikzlibrary{arrows,positioning} 
\makeatother

\usepackage{babel}
  \providecommand{\examplename}{Example}

\begin{document}
\begin{center}
\textbf{\Large{An overview of type theories}}
\par\end{center}{\Large \par}

\begin{center}
Nino Guallart
\par\end{center}

\title{}
\begin{abstract}
Pure type systems arise as a generalisation of simply typed lambda
calculus. The contemporary development of mathematics has renewed
the interest in type theories, as they are not just the object of
mere historical research, but have an active role in the development
of computational science and core mathematics. It is worth exploring
some of them in depth, particularly predicative Martin-Löf's intuitionistic
type theory and impredicative Coquand's calculus of constructions.
The logical and philosophical differences and similarities between
them will be studied, showing the relationship between these type
theories and other fields of logic. 

\end{abstract}
\maketitle

\section{Introduction}

This paper is an overview of generalised type systems, in particular
normalising dependent systems, focusing on a comparison between predicative
and impredicative dependent theories. It is intended as a very basic
introduction, so no previous knowledge of the topic is assumed. The
first part of the paper is dedicated to preliminary notions of $\lambda-$calculus
and simple type theory, and the second one to the questions regarding
generalised type systems.

Type theories arose as an alternative to set theory due to some contradictory
situations that appear in naïve set theory when considering certain
definitions of sets. The most important example is Russell's paradox,
which involves the set of all sets that are not members of themselves,
$R=\{x|x\notin x\}$ (Cfr. (Russell, 1903) and (Russell, 1980)). It
is clear that this set is a member of itself if and only if it is
not a member of itself, which is contradictory. The discovery of this
paradox revealed a fatal flaw in naïve logic and forced to reconsider
their basic principles. Logicians and mathematicians tackled this
problem by restricting the way sets can be formed with two basic approaches:
\begin{itemize}
\item Formulation of axiomatised set theories (ZFC, NBG...). Instead of
relying on a principle of unrestricted comprehension, these systems
have a version of an axiom or rule of separation, which needs a previously
defined set in order to build a new one from it and a certain predicate
ranging over its elements.
\item Type theories. In layman's terms, a type is almost the same thing
as a set, except that types form a hierarchy that avoids self-reference,
since a type contains elements of a lower range. In this paper we
will cover some aspects of the development of type theories.
\end{itemize}
Self-reference, which plays a crucial part in many paradoxes like
this, is closely related to impredicativity; a definition is impredicative
if it quantifies over a set containing the entity being defined, and
predicative otherwise. Russell's paradox is impredicative, and some
authors like Russell himself or Poincaré thought that impredicativity
was problematic, arguing that impredicative definitions lead to a
vicious circle. However, many basic mathematical definitions are impredicative
(e.g. least upper bound) and Ramsey (1931) stated that many impredicative
definitions are actually harmless and non-circular. Predicative mathematics
avoids any problem regarding impredicativity by relying only on predicative
principles, but the rejection of impredicative definitions leaves
out some very important ones, such as power sets, making thus the
task of predicative mathematics quite difficult. The the vast majority
of mathematicians have continued being impredicative, but we will
see later how new predicative theories emerge again in generalised
type theories.

Although several authors such as Russell, Ramsey or Gödel created
their own versions of type theory, probably the most remarkable one
is Church's system. Its main feature is that it is formulated in $\lambda$-calculus.
From now on we will focus on type theories based on it, that is, typed
$\lambda$-calculus. In order to understand them properly, we will
study some basic notions of untyped $\lambda$-calculus, before dealing
with typed theories.

\section{Simple type theory }

\subsection{Untyped lambda calculus }

Lambda calculus ($\lambda$C) is a formal system created by Church
in order to deal with the notions of recursion and computability in
connection to the problem of halting (\emph{Entscheidungsproblem}).
Instead of relying on elements and sets as primitives, it uses $\lambda$-terms,
which are used to explore the concept of function.

The set of $\lambda$-terms, or just terms in short, is defined recursively
from a countable set of variables $Var$. Using Backus-Naur notation,
formation of terms can be expressed in the following scheme, where
$x$ is any variable, $E$ and $F$ are arbitrary terms and the dots
denote optional constant terms:

\[
E,F:=x|\lambda x.E|(EF)|\ldots
\]

In other words, the set $\Lambda$ of terms can be defined recursively
as follows:
\begin{itemize}
\item If $x\in Var$, $x\in\Lambda$. 
\item If $E\in\Lambda$, $F\in\Lambda$, then $(EF)\in\Lambda$. 
\item If $E\in\Lambda$ and $x\in Var$, $\lambda x.E\in\Lambda$.
\end{itemize}
Simplification rules for parentheses can be applied as usual: outer
parentheses are omitted and left association is assumed, so $EFG$
stands for $((EF)G)$; abstraction is right associative, so $\lambda x.\lambda y.E$
stands for $\lambda x.(\lambda y.E)$; abstractions can be expressed
in contracted mode, for example $\lambda xy.E$ instead of $\lambda x.\lambda y.E$. 

Variables may occur bound in a term by the $\lambda$ abstraction
or be free: in $\lambda x.xy$ $x$ is bound and $y$ is free. A basic
and important opperation is the substitution of variables for terms.
We will write $P[x:=t]$ to indicate the substitution of the variable
$x$ in $P$ for the term $t$. For example, $xyx[x:=z]$ and $xyx[x:=ab]$
gives out $zyz$ and $abyab$ respectively, since all $x$ have been
substituted by $z$ and $ab$. Substitution is the basis of $\alpha-$conversion
and $\beta$-reduction.

\subsubsection*{$\alpha$-conversion.}

Variable names are conventional, so any variable in a term can be
renamed, that is, substituted by another variable which does not occur
in the term. This is called $\alpha-$conversion: $\lambda x.x$ can
be converted to $\lambda y.y$, so they are equivalent under $\alpha-$conversion,
$\lambda x.x\rightarrow_{\alpha}\lambda y.y$. Since the new variables
must not occur in the term, $\lambda x.xy\rightarrow_{\alpha}\lambda y.yy$
is wrong, whereas $\lambda x.xy\rightarrow_{\alpha}\lambda z.zy$
is correct.

\subsubsection*{$\beta$-reduction and operational semantics of $\lambda$C}

Term evaluation in $\lambda$C is based on $\beta$-reduction, the
application of functions over their arguments substituting the bound
variables in the function by their corresponding arguments (actual
process may be more complex, but for the sake of simplicity this work
will focus solely on $\beta$-reduction as the basis of term evaluation):

\[
(\lambda x.E)F\Longrightarrow E[x:=F]
\]

\[
(\lambda x.xy)a\Longrightarrow xy[x:=a]\Longrightarrow ay
\]

After a $\beta$-reduction, a term may or may not need ulterior reduction.
A reducible subterm within a term is usually called a reducible expression
or a redex in short. If a term cannot be reduced further, it is said
to be in its normal form, or to be a value. In $\lambda$-calculus,
any term evaluation converges to a value or diverges. Untyped lambda
calculus is not normalising, that is, it is not always possible to
reach a value. 
\begin{example*}
Identity function $I:=\lambda x.x$, when applied to a term, gives
the same term as an output:

\[
(\lambda x.x)e\Longrightarrow x[x:=e]\Longrightarrow e
\]

Since $I$ is also a term, it can be applied to itself, giving itself
as an output:

\[
II\Longrightarrow I
\]

\end{example*}
\[
(\lambda x.x)(\lambda x.x)\Longrightarrow x[x:=\lambda x.x]\Longrightarrow\lambda x.x
\]

\begin{example*}
$\omega$ operator $\lambda x.xx$ takes an input $e$ and gives $ee$
as an output; this new term may or may not be reduced again. If applied
to itself, the result is the same as the original term and therefore
a value is never reached.
\end{example*}
\[
\omega\omega\Longrightarrow\omega\omega
\]

\[
(\lambda x.xx)(\lambda x.xx)\Longrightarrow xx[x:=\lambda x.xx]\Longrightarrow(\lambda x.xx)(\lambda x.xx)
\]

\begin{example*}
Sometimes, the reduction of a term gives a more complex term and thus
a value is never reached:
\end{example*}
\[
(\lambda x.xxx)(\lambda x.xxx)\Longrightarrow(\lambda x.xxx)(\lambda x.xxx)(\lambda x.xxx)\Longrightarrow
\]

\[
\Longrightarrow(\lambda x.xxx)(\lambda x.xxx)(\lambda x.xxx)(\lambda x.xxx)\Longrightarrow\ldots
\]

Different strategies can be applied to reduce a term. The normal way
of reducing terms is named \emph{call-by-value}, which means reducing
the rightmost redex first. The opposite strategy, reducing the leftmost
redex first, is named \emph{call-by-name}. There are many other strategies.
However, if the term eventually reaches a normal form, it holds the
Church-Rosser property, that is, the normal form is independent of
the strategy that has been chosen. It is worthwhile noting that this
semantics of $\lambda$C is operational, since it describes the way
$\lambda$-terms function, not what they denote.

Recursion in $\lambda$C is possible using fixed point combinators,
which are terms $R$ such as, when over another term $T$, give

\[
RT\Longrightarrow T(RT)\Longrightarrow T(T(RT))\Longrightarrow\ldots
\]

If we make $Q:=RT$ then $Q:=TQ$, which is a fixed point, hence its
name. There are infinite fixed point combinators, although the most
known is Curry's $Y$ combinator.

The following closing remarks can be made: 
\begin{itemize}
\item Lambda calculus' foundations rely on the use of abstraction and application:
$\lambda x.P$ is an anonymous function which takes an input $x$
and gives $P(x)$ as an output. The application $PQ$ is the result
of considering term $Q$ as the input of $P$.
\item Recursion allows the construction of complex formulae and is the basis
of the computational power of lambda calculus.
\item $\lambda$C is equivalent to a Turing machine, and therefore is an
abstract model of a programming language (Landin, 1965). In $\lambda$C,
Boolean values are possible and natural numbers can be defined in
the form of Church's numerals, and there are also logical and arithmetical
operations over them and flux control operators. It is not decidable,
so in principle it is not possible to know whether the evaluation
of a term will end or not. Whereas TM is the model of imperative programming
languages, $\lambda$C is the basis of functional ones. 
\end{itemize}

\subsection{Simply typed lambda calculus }

Simply typed lambda calculus was also originally developed by Church
(1940,1941). It is a higher order logic system based on lambda calculus
and it uses the same syntax. We will not see the original formulation
of Church here, but a more recent one which will allow us to connect
with current developments.

A simply typed lambda calculus (STLC) has a non-empty set of base
types. The rest of the types, which are function types, can be built
from them by the application of the type constructor $\rightarrow$.
Types and terms can be defined by recursive rules, where $E$ and
F are any terms, $\sigma$ and $\tau$ any types, $x$ any variable,
$\beta$ a base type, and the dots denote optional terms belonging
to a given type:

\[
\sigma,\tau:=\beta|\sigma\rightarrow\tau
\]

\[
E,F:=x|\lambda x:\sigma.E|EF|\ldots
\]

$\lambda x:\sigma.E$ is also written $\lambda x^{\sigma}.E$, the
first form will be preferred, although the second one will be used
sometimes in order to increase clarity. There are two styles of typing,
Church style and Curry style. Although in a subtle way, their differences
are not only syntactical, but also semantical. In Curry style, variable
types are not explicitly stated, whereas in Church style it is imperative
to declare the type of every variable. The problems with Curry style
will not be discussed here, and this entire section refers only to
Church style. The following paragraphs will focus on the difference
between well-formed terms (actually they are called \emph{pre-terms})
and well-typed terms. 

In typed theories there are \emph{typing environments}, also called
variable assignments, which are (possibly empty) sets of associations
between types and variables, so each variable has its type, and we
write $x:\sigma$. In a typing environment, a \emph{typing judgment}
$\Gamma\vdash E:\sigma$ is a statement of the fact that, under environment
$\Gamma$, term $E$ is well-formed and has type $\sigma$. In other
words, it asserts that $E$ is a term of type $\sigma$ if and only
if its free variables are of the types specified in the typing environment.

\subsubsection*{Typing rules.}

The associations between terms and types are maintained by \emph{typing
rules}, which are inference rules from a group of premises to a conclusion,
all of them being typing judgments. The following typing rules specify
how to assign a type to a certain well-formed syntactic construction,
where \noun{var} stands for variable (a variable within a context
is a well-typed term), \noun{abs} for abstraction (given a variable
and a well-typed term, the $\lambda-$abstraction of the term is well-typed)
and \noun{app} for application (the application of two well-typed
terms is also well-typed): 

\medskip{}

\noindent %
\begin{minipage}[t]{1\columnwidth}%

\AxiomC{$x:\sigma \in\Gamma $} 
\RightLabel{\scriptsize VAR}
\UnaryInfC{$\Gamma \vdash x:\sigma $} 
\DisplayProof\hfill{}
\AxiomC{$\Gamma , x: \sigma \vdash e: \tau$} 
\RightLabel{\scriptsize ABS}
\UnaryInfC{$\Gamma \vdash (\lambda x:\sigma .e ): (\sigma \rightarrow \tau)$} 
\DisplayProof%
\end{minipage}

\begin{prooftree} 
\AxiomC{$\Gamma \vdash e: \sigma \rightarrow \tau$} 
\AxiomC{$\Gamma \vdash f: \sigma $} 
\RightLabel{\scriptsize APP}
\BinaryInfC{$\Gamma \vdash ef: \tau $} 
\end{prooftree}

In STLC well formed terms may be not well-typed, and the typability
of terms are given by the previous rules, so the set of well-typed
terms is a proper subset of the set of all well formed terms. For
example, $(\lambda x:\sigma.x)(y:\tau)$ is well formed since it can
be built following the mentioned syntactic rules of formation, but
not well-typed, because it does not satisfy rule APP, since $\lambda x:\sigma.x$
has type $\sigma\rightarrow\sigma$ and $y$ has type $\tau$. In
this case, its evaluation gives as an output an error. 

Provided that the types are right, term evaluation is based on $\beta$-reduction,
like in untyped $\lambda-$calculus:

\medskip{}

\begin{prooftree} 
\AxiomC{$(\lambda x:\sigma.E)F$} 
\RightLabel{\scriptsize E-APP}
\UnaryInfC{$E[x:=F]$} 
\end{prooftree}

\noindent %
\begin{minipage}[t]{1\columnwidth}%

\AxiomC{$t \rightarrow t'$} 
\RightLabel{\scriptsize E-APP1}
\UnaryInfC{$tu \rightarrow t'u$} 
\DisplayProof\hfill{}
\AxiomC{$t \rightarrow t'$} 
\RightLabel{\scriptsize E-APP2}
\UnaryInfC{$ut \rightarrow ut'$} 
\DisplayProof%
\end{minipage}

\medskip{}

We can also consider the void type 0 with no terms and the unit type
with a single term, $*:1$ such as for any type $T$, it holds $\lambda x^{T}.*:1$.
Since type 0 has no terms, it is also easy to see that we can not
construct a valid term of type $T\rightarrow0$.

STLC can be extended to STLC with pairs, in which the product type
$\sigma\times\tau$ appears, which is the type of pairs of terms $(s,t)$.
The notation for types and terms in a system with two type constructors,
$\rightarrow$ and $\times$ is the following:

\[
\sigma,\tau:=\beta|\sigma\rightarrow\tau|\sigma\times\tau
\]

\[
E,F:=x|\lambda x:E\sigma.F|EF|(E,F)
\]

The rules for the formation and elimination of pairs are the following:

\begin{prooftree} 
\AxiomC{$\Gamma \vdash e:\sigma $} 
\AxiomC{$\Gamma \vdash f:\tau $} 
\RightLabel{\scriptsize I-PAIR}
\BinaryInfC{$\Gamma \vdash (e,f):\sigma \times \tau$} 
\end{prooftree}

\noindent %
\begin{minipage}[t]{1\columnwidth}%

\AxiomC{$\Gamma \vdash (e,f):\sigma \times \tau $} 
\RightLabel{\scriptsize E-PAIR L}
\UnaryInfC{$\Gamma \vdash \pi_{1} (e,f):\sigma $} 
\DisplayProof\hfill{}\AxiomC{$\Gamma \vdash p:\sigma \times \tau $} 
\RightLabel{\scriptsize E-PAIR R}
\UnaryInfC{$\Gamma \vdash \pi_{2} (e,f):\tau $} 
\DisplayProof%
\end{minipage}

\medskip{}

$\pi_{1}$ and $\pi_{2}$ are respectively the first and second projection
of the pair, so $\pi_{1}(e,f)=e$ and $\pi_{2}(e,f)=f$. Both STLC
and STLC with pairs are strongly normalising, that is, the evaluation
of a term eventually gives a value if it is well-typed or an error
if not, and therefore they are decidable. In this way, STLC and its
derivatives may serve as a basis for type checking, which has its
practical, computational applications. Not all well-formed terms are
well-typed, unlike untyped $\lambda$C, which has no typing restrictions.
Because of it, STLC and STLC with pairs are not Turing-complete, since
they are not equivalent to $\lambda$C, and they are less expressive
than it.

\subsubsection*{Type theories and category theory. }

Until now, we have considered syntactical and operational aspects
of type theories, neglecting their semantical facet. Now we are going
to see the semantical application of category theory to type theories.
Category theory was developed by Eilenberg and MacLane (1945) for
the study of algebraic topology, but it soon showed to be useful in
many other fields. Its connection to logic has been studied since
the the works of Lawvere (1963); in this paper we will also follow
the works of Lambek and Scott (1986).Lawvere and other authors have
also studied category theory as a new approach to the foundations
of mathematics, dealing with abstract mathematical structures, instead
of using the traditional set-theoretical frame. 

A category $\mathcal{C}$ is made of objects and morphisms or arrows
between objects. Each morphism has as a domain a given object and
as codomain another one, so for example morphism$f:A\rightarrow B$
has $dom(f)=A$ and $cod(f)=B$. The composition of arrows is associative,
$(fg)h=f(gh)$, and for every object $X$ there is an identity arrow
$1_{X}:X\rightarrow X$ such as for an arbitrary pair of morphisms
$g:Y\rightarrow X$ and $f:X\rightarrow Z$, $f1_{X}=f$ and $1_{X}g=g$.
There can be additional compositions of objects that will be useful
later, such as the product $X\times Y$ of two objects $X$ and $Y$,
which is an object $Z=X\times Y$ with two projection morphisms $\pi_{1}$
and $\pi_{2}$ satisfying $\pi_{1}Z=X$ and $\pi_{2}Z=Y$; the exponential
or power object $X^{Y}$ can be seen as the class of all morphisms
from $Y$ to $X$ (the actual definition is more complex, but for
the purposes of this paper, this definition is valid). A category
has a given structure that can be simple or complex, and we can be
consider categories made of categories. We call \emph{functors} the
mappings between these categories preserving the structure from one
category into the other.

There is a bidirectional relationship between categories and type
theories (Cfr. Asperti and Longo, 1991). Type theories can be interpreted
using category theory, and conversely, we can formalise categories
in the language of type theories. Generally, we can say that the relation
between a type theory and its corresponding category is akin to \emph{syntax
vs semantics}. We can show this studying the relationship between
STLC with pairs and Cartesian closed categories (Cfr. Lambek and Scott,
1988). 

A category is a Cartesian closed category (CCC) if and only if it
has a terminal object $T$ such as for every object $X$ there is
a unique morphism $X\rightarrow T$, and if for any two objects $X$
and $Y$ there exists the product $X\times Y$ and the exponential
object $X^{Y}$. We can see that the elements of STLC $\mathcal{T}$
form a CCC: we can consider base types of $\mathcal{T}$ as objects,
and it is easy to see that type 1 can be interpreted as the terminal
object $T$. Constant terms $a$ of a given type $A$ can be seen
as morphisms $a:1\rightarrow A$ from the terminal object to their
corresponding object $A$, and more broadly, terms can be interpreted
as morphisms, whereas application of terms equates to composition
of morphisms. For every object $A$ there is an identity function
$1_{A}$ which can be seen as the identity morphism for $A$, and
it can be proved that composition of terms is associative. Product
type $A\times B$ equates to the product of objects, and function
type $A\rightarrow B$ is interpreted as the power object $B^{A}$,
which is the class of equivalence of all morphisms from $A$ to $B$
with a free variable of type $A$. Therefore, we have that the interpretation
of $\mathcal{T}$, $Syn(\mathcal{T})$, is a category, and more precisely
a CCC which we can call $\mathcal{C}'$.

Conversely, we can see that the structure of a given CCC $\mathcal{C}$
can be described using the language of a STLC $\mathcal{T}'=Lang(\mathcal{C})$,
where $Lang(\mathcal{C})$ is the smallest type theory that preserves
the structure of the category, that is, $\mathcal{T}'$ is the internal
language of $\mathcal{C}$ (Cfr. Johnstone (2003)). We can denote
the terminal object $T$ with the unit type, and assign a type to
each object in the category. Product of objects and exponential objects
equate to product types and function types respectively, and for each
morphism $f$ between objects $X$ and $Y$, there is a term $\lambda x^{X}.f$
with a type $X\rightarrow Y$.

Actually, this relation between categories and type theories is subtler,
since given a CCC $\mathcal{C}$ and a STLC $\mathcal{T}$ with the
same structure, the internal language of the category $Lang(C)$ is
not equal to $\mathcal{T}$, but homomorphical to it ($Lang(\mathcal{C})\cong\mathcal{T}$),
and the same with the semantical interpretation of the theory, $Syn(\mathcal{T})\cong\mathcal{C}$.
$Syn$ and $Lang$ are actually two functors between categories, so
we can succintly write their relationship in this way: $\mathcal{C}\overset{Syn}{\underset{Lang}{\leftrightharpoons}}\mathcal{T}$. 

\begin{center}
\begin{tabular}{|c|c|}
\hline 
\noun{Type theory} & \noun{Category theory}\tabularnewline
\hline 
\hline 
{\small{Types}} & {\small{Objects}}\tabularnewline
\hline 
{\small{Unit type ($\top$ or 1)}} & {\small{Terminal object}}\tabularnewline
\hline 
{\small{Product type $A\times B$}} & {\small{Product of objects $A\times B$}}\tabularnewline
\hline 
{\small{Function type $A\rightarrow B$}} & {\small{Exponential object $B^{A}$}}\tabularnewline
\hline 
{\small{Terms}} & {\small{Morphisms}}\tabularnewline
\hline 
{\small{Pair of terms $(f,g)$}} & {\small{Pair of morphisms $(f,g)$}}\tabularnewline
\hline 
{\small{Projections of terms, $\pi_{1}$and $\text{\ensuremath{\pi}}_{2}$}} & {\small{Projections of morphisms, $\pi_{1}$and $\text{\ensuremath{\pi}}_{2}$}}\tabularnewline
\hline 
{\small{Abstraction $\lambda x^{A}.f:B$}} & {\small{Arrow $f:A\rightarrow B$ with a free variable $x:A$}}\tabularnewline
\hline 
{\small{Application $fg$}} & {\small{Composition of arrows $fg$}}\tabularnewline
\hline 
\end{tabular}
\par\end{center}

\begin{center}
\noun{Table 1.}
\par\end{center}

\medskip{}

\section{The lambda cube and generalised type systems }

From the 70s onwards there have been prominent new works in the field
of type theories. A complex family of typed lambda systems, each with
their own features, flaws and strengths, has been developed from Church's
simple type theory. There are several reasons that explain this renewed
interest. Firstly, some practical reasons: type theory is closely
related to proof theory and therefore to computing and typechecking.
Secondly, from a theoretical point of view, the rise of category theory
and further developments into their connections to logic have been
the object of intensive research, as type theories can be studied
from a categorical point of view. In relation to this, there have
been renewed efforts in the research of type theories as an alternative
foundation of mathematics, particularly in constructive mathematics,
mainly since the works of Per Martin-Löf, and new fields within logic
and mathematics such as homotopical type theory have appeared, advancing
this study.

\subsubsection*{Pure type systems. }

STLC can be considered as the basis of a family of type systems that
have been named pure type systems (PTS) or generalized type systems
(GTS). They are a group of type systems that, unlike STLC, allow
dependencies between types and terms. Broadly speaking, the main difference
between a simply typed theory and a pure type system is that the latter
allows judgments over types. In STLC types and terms are two disjoint
groups, in PTS this distinction is blurred or erased. 

Instead of creating confusing categories such as types of types, types
of types of types and so on, the concepts of kind and sort are used
instead, which are generalisations of the notion of type. $\star$,
also called \emph{Prop}, is any usual type of terms or propositions.
$\square$ or \emph{Type} denotes the higher-order type of a type.
The set of sorts is thus $S=\{\star,\square\}$. A kind $\star\rightarrow\square$
maps from a type of terms to a certain kind of types and so on. In
a given context, it can be asserted that $\sigma:\square$ for specifying
that $\sigma$ is a valid type in the same way that $e:\sigma$ .
The scheme for kinds and expressions in PTS is the following, where
$V$ is any variable and $S$ any sort and $K$ any kind ($\Pi$ operator
will be explained later):

\[
K:=\star|\square|K\rightarrow K
\]

\[
T,U:=V|S|TU|\lambda V:T.U|\Pi V:T.U
\]

Not all PTS are normalising, that is, when evaluated not all of them
reach a value. In the next subsections we will study a more reduced
family of generalised type theories, all of which are decidable and
normalising.

\subsection{Barendregt's lambda cube}

Barendregt (1991) considers four possible relations between terms
and types: terms depending on terms, types on types, terms on types
and types on terms. If we omit the first one, we have three possibilities
that can be represented as axes of a cube, other features such as
subtyping could be represented in additional dimensions (also Cfr.
(Baredregt, 1992)): 

\begin{center}
\begin{tikzpicture}
\node (a) at (0,0) {$\lambda_{\rightarrow}$};
\node (b) at (3,0) {$\lambda_{P}$};
\node (c) at (0,3) {$\lambda_{2}$};
\node (d) at (3,3) {$\lambda_{P2}$};

\node (a2) at (1.7,1.7) {$\lambda_{\underline\omega}$};
\node (b2) at (4.7,1.7) {$\lambda_{P\underline\omega}$};
\node (c2) at (1.7,4.7) {$\lambda_{\omega}$};
\node (d2) at (4.7,4.7) {$\lambda_{P\omega}$};

\draw [thick,->] (a) -- (b);
\draw [thick,->] (a) -- (c);
\draw [thick,->] (b) -- (d);
\draw [thick,->] (c) -- (d);

\draw [->] (a2) -- (b2);
\draw [->] (a2) -- (c2);
\draw [thick,->] (b2) -- (d2);
\draw [thick,->] (c2) -- (d2);

\draw [->] (a) -- (a2);
\draw [thick,->] (b) -- (b2);
\draw [thick,->] (c) -- (c2);
\draw [thick,->] (d) -- (d2);

\end{tikzpicture}
\par\end{center}

\begin{center}
\noun{Figure 1.}
\par\end{center}

These three possibilities considered by Barendregt are the following:
\begin{itemize}
\item Type polymorphism (terms depending on types, $\square\rightarrow\star$):
universal quantification over types in order to be able to define
type variables. System $\lambda$2, also called system F, which was
discovered independently by Girard (1972) and Reynolds (1974), is
a second order lambda calculus or polymorphic lambda calculus.\end{itemize}
\begin{example*}
In simple type theory there is no unified identity function, unlike
$I:=\lambda x.x$ of untyped $\lambda$C. Instead of a single function,
in STLC there is a family of identity functions, one function $\lambda x:\sigma.x$
for each type $\sigma$. If quantification over types is allowed,
then the previous idea can be formalised succintly in the way system
F does:
\end{example*}
\[
\Lambda\alpha.\lambda x^{\alpha}.x:\ \forall\alpha.\alpha\rightarrow\alpha
\]

Analogously to the way simply typed lambda calculus defines types
and terms, so does system F, being $\beta$ a base type and $\alpha$
any type variable:

\[
\sigma,\tau:=\beta|\alpha|\sigma\rightarrow\tau|\forall\alpha.\sigma
\]

\[
E,F:=x|EF|\lambda x:\sigma.E|\Lambda\alpha.E
\]

Unlike in STLC, types can appear in terms, like the case of the previous
example $\Lambda\alpha.\lambda x^{\alpha}.x$. 

Two new rules for introduction ($\forall I$) and elimination ($\forall E$)
of generalisation over types:

\begin{prooftree} 
\AxiomC{$\Gamma \vdash M:\tau $} 
\AxiomC{$\alpha \not \in \Gamma$} 
\RightLabel{\scriptsize $\forall$ I}
\BinaryInfC{$\Gamma \vdash \Lambda \alpha .M : \forall \alpha .\tau $} 
\end{prooftree}

\begin{prooftree} 
\AxiomC{$\Gamma \vdash M: \forall \alpha .\sigma $} 
\AxiomC{$\Gamma \vdash T:\tau$}
\RightLabel{\scriptsize $\forall$ E}
\BinaryInfC{$\Gamma \vdash M T : \sigma [\alpha:=\tau]$} 
\end{prooftree}
\begin{itemize}
\item Type constructors (types depending on types, $\square\rightarrow\square$):
abstraction of new types from previous ones. This is a remarkable
implementation, because new types are built within the language, not
the metalanguage of kinds (Roorda, 2000). There will be rules for
kinds and expressions (which comprise both types and terms) instead
of for types and terms. $\lambda_{\omega}$ calculus (Girard, 1972)
and system $F_{\omega}$/$\lambda2_{\omega}$ use type constructors.
\end{itemize}
Rules for this implementation can be quite complex and they will not
be covered here. Product type is a simple example of this feature.
Other type constructors such as list constructors or higher order
types are within this category.
\begin{itemize}
\item Dependent types (types depending on terms, $\star\rightarrow\square$):
The last possibility is building types depending on previous types.
Here are two possibilities, product types and dependent sum (or pair)
types. Product types ($\Pi$-types) generalise the idea of universal
quantification. Considering a non-empty type $A$ that will serve
as an index, $\Pi$-types generate a family of types $B(a)$ depending
on every $a\in A$. If $B(x)$ is a constant type $B$, then $\Pi_{x:A}B(x)$
equals $A\rightarrow B$. Conversely, sum types ($\Sigma$-types)
in the form $\Sigma_{x:A}B(x)$, are the type of pairs of terms of
the form $(a^{A},b^{B[x:=A]})$ in which the type of the second element
depends on the value of the first term. Sum types are the equivalent
of existential quantification and, if $B$ is a constant type, $\Sigma_{x:A}B(x)$
equals $A\times B$. Two examples of systems with dependent types
are system LF, which is STLC with dependent types, and calculus of
constructions, that will be covered soon.
\end{itemize}
To sum up:

\begin{center}
\begin{tabular}{|c|c|c|}
\hline 
System & Relations & Examples\tabularnewline
\hline 
$\lambda_{\rightarrow}$ & $\star\rightarrow\star$ & STLC\tabularnewline
\hline 
$\lambda_{2}$ & $\star\rightarrow\star$,$\square\rightarrow\star$ & System F\tabularnewline
\hline 
$\lambda\underline{\omega}$ & $\star\rightarrow\star$,$\square\rightarrow\square$ & Weak $\lambda_{\omega}$\tabularnewline
\hline 
$\lambda_{\omega}$ & $\star\rightarrow\star$,$\square\rightarrow\star$,$\square\rightarrow\square$ & System F$\omega$\tabularnewline
\hline 
$\lambda P$ & $\star\rightarrow\star$,$\star\rightarrow\square$ & System LF\tabularnewline
\hline 
$\lambda P_{2}$ & $\star\rightarrow\star$,$\star\rightarrow\square$,$\square\rightarrow\star$ & $\lambda P_{2}$\tabularnewline
\hline 
$\lambda P_{\underline{\omega}}$ & $\star\rightarrow\star$,$\star\rightarrow\square$, $\square\rightarrow\square$ & Weak $\lambda P_{\omega}$\tabularnewline
\hline 
$\lambda P_{\omega}$ & $\star\rightarrow\star$, $\square\rightarrow\star$, $\star\rightarrow\square$,
$\square\rightarrow\square$ & CoC\tabularnewline
\hline 
\end{tabular}
\par\end{center}

\begin{center}
\noun{Table 2.}
\par\end{center}

\medskip{}

All application and abstraction rules of these systems can be summarised
as follows if we consider the whole cube system (Cfr. (Roorda, 2000)):

\begin{prooftree} 
\AxiomC{$\Gamma, x:A \vdash b:B$} 
\AxiomC{$\Gamma \vdash (\Pi x:A.B) :t \in \{\star,\square\}$} 
\RightLabel{\scriptsize ABS}
\BinaryInfC{$\Gamma \vdash (\lambda x:A.b) : \Pi x:A.B$} 
\end{prooftree}

\begin{prooftree} 
\AxiomC{$\Gamma \vdash f: (\Pi x:A.B) :t$} 
\AxiomC{$\Gamma \vdash a:A $} 
\RightLabel{\scriptsize APP}
\BinaryInfC{$\Gamma \vdash fa : B [x:=A]$} 
\end{prooftree}

Calculus of constructions, which will be mentioned later, allows all
of these extensions and, therefore, it is placed in the uppermost
right back the cube. Other systems such as system F only allow some
of them, being placed in other corners. All of them are strongly normalising,
that is, end giving a value or an error in a finite number of steps.

Semantically, $2$-categories (categories over morphisms) or more
generally $n$-categories and higher-order type theories are related
in a way similar to the link between CCC and STLC. Semantics of these
type theories is far more complicated than the semantics of STLC and,
therefore, just a brief sketch will be made. The starting point of
the categorical interpretation of dependent types is the consideration
of slice categories, which are categories in which the objects are
morphisms over a given object. It is easy to see the relationship
within slice categories and dependent types, since an object $A$
in which every morphism $a$ gives out a slice category $B(a)$ can
be considered the semantical interpretation of the product type $\Pi x:A.B$.
\emph{Locally closed Cartesian categories} (lCCC) are categories in
which all slice categories are CCC and a kind of dependent type theory
is the internal language of lCCC (Seely, 1984), in the same way that
STLC with pairs is the internal language of CCC.

\subsection{Calculus of constructions}

Coquand's calculus of constructions (CoC in short) is a higher-order
lambda calculus theory that combines polimorphism and type construction
of Girard's system $F_{\omega}$ with dependent types. The syntax
of its kinds and terms is the following: 

\[
K:=\star|\square_{i}\ (i\geq1)
\]

\[
\sigma,\tau,M,N:=x|K|\Pi x^{\sigma}.\tau|\lambda x^{\sigma}.M^{\tau}|MN{}_{\tau[x:\sigma]}
\]

CoC distinguishes between the impredicative type of predicates ($\star$,
small types), a predicative hierarchy of types of types ($\square_{i}$,
large types), and the type of all large types. As in other PTS, from
$\Pi$ constructor it is possible to derive the usual logical operators.
CoC also has void type and $1$ type, and it is easy to create types
for truth values and natural numbers. It has several variations, such
as CoC with inductive types, but they will not be treated here.

So far, all considered type theories are impredicative. None of them
are problematic, since they are strongly normalising. For example,
the polymophic identity in System F $\Lambda\alpha\lambda x:\alpha.x$
can take as arguments its own type $\forall\alpha.\rightarrow\alpha$
and then itself, but in a way that leads to no circularity. Some impredicative
theories such as Girard's system U are inconsistent, but this is not
the case of the type systems considered in the $\lambda-$cube.

\subsection{Intuitionistic type theory}

Per Martin-Löf's intuitionistic type theory (ITT in short) allows
to introduce contemporary predicative theory types in this discussion.
It can be considered and extension of STLC with higher order predicates
and quantification over types within a mathematical constructivist
programme. Strictly speaking, ITT does not belong to the $\lambda$-cube,
but it has an expressive power similar to the one of CoC, so a brief
comparison between them seems reasonable. Semantically, Seely (1984)
showed that there is a relationship between locally Cartesian closed
categories and ITT. It has shown some prominent features in the field
of programming due to its connections to proof theory, and it also
aims to serve as a constructivist theory for the foundations of mathematics.
Homotopy type theory expects to follow this aspiration, since the
programme of the Univalent Foundations of Mathematics conceive this
field as an extension of ITT with a homotopical interpretation (Voevodsky,
2013). 

As has been said before, previous type theories are impredicative,
yet this feature is unproblematic. ITT's first formulation (Martin-Löf,
1971) was also impredicative, but was soon discovered to be inconsistent.
Later developments (Martin-Löf, 1975) avoided the problems of this
earlier version with a predicative formulation, a common feature in
other constructivist approaches to the foundations of mathematics.
Categorically, predicative theories have a more general structure;
whereas many impredicative theories lie on topos theories, predicative
ones rely on structures such as pretopos, which are more broad categories
but require more complex proofs, since impredicative definitions are
rejected.

The concept of universe in ITT, which is similar to Grothendieck's
universe, is crucial in this theory, which is in consonance with the
constructive theses of its author. It can be considered as a closed
type of types built according to certain conditions. An earlier version
of ITT had an impredicative universe, but it showed to be inconsistent.
Later versions (Cfr. Martin-Löf, 1975) consider a hierarchy of predicative
universes $(U_{i},Type_{i})$ or $(U_{i},T_{i})$ in which a given
$Type_{i}$ has a type $Type_{i+1}$. 

\begin{prooftree} 
\AxiomC{$ x \in U_i$} 
\RightLabel{\scriptsize U-FORM}
\UnaryInfC{$x: Type_i$} 
\end{prooftree}

\begin{prooftree} 
\AxiomC{$ x \in U_i$} 

\RightLabel{\scriptsize U-FORM}
\UnaryInfC{$Type_i: Type_{i+1}$} 
\end{prooftree}

Like CoC, ITT uses dependent types, and universes are closed under
operations. Intuitionistic type theory can be formalised stating its
own typing context and its typing rules and judgments, in the same
way as other type theories. There are also several predefined finite
types, void type (0), unity or truth (1), and bool (2). Both CoC and
ITT are strongly normalising and therefore non Turing complete.

\subsection{Curry-Howard isomorphism}

From the works of Curry (1934) and Howard (1969) on, it has been established
the correspondence between each type theory and a style of logical
calculus, or more broadly between type theory and proof theory. In
simply typed lambda calculus, types can be built in a way akin to
predicate logic well-formed formulae. This establishes a link between
type theories and logical calculi, since the types of different systems
can be treated as well-formed formulae of the corresponding logical
systems. We will see this covering mainly the relationship between
STLC and intuitionistic predicate logic.

In intuitionistic logic, the semantics is given by the Brouwer-Heyting-Kolmogorov
interpretation: truth is identified with provability, so saying that
$A$ is true means that there is a proof $a$ for $A$. If we interpret
types as predicates and terms as proofs, we can see that these two
judgments are equivalent:
\begin{itemize}
\item Proof $a$ proves the predicate $A$.
\item Type $\alpha$ can be built and this type is inhabited by a term $a:\alpha$. 
\end{itemize}
Let $0$ be the void type (the type with no proofs/terms) and $\neg\alpha\equiv_{def}\alpha\rightarrow0$,
that is, $\alpha$ leads to contradiction; $\alpha\times\beta$ means
having a pair of terms, $a:\alpha$ and $b:\beta$, so it can be identified
with $\alpha\land\beta$; $\alpha+\beta$ is having either a term
$a:\alpha$ or a term $b:\beta$, that is, $\alpha\lor\beta$. The
rule of term application ($e:\sigma\rightarrow\tau,f:\sigma\vdash ef:\tau$)
equates to \emph{modus ponens}. We can see that there is an isomorphism
between STLC and predicate logic. However, the constructed predicate
logic is not classical predicate logic, but an intuitionistic system,
and some classical tautologies such as Peirce's law cannot be obtained
unless extra axioms are added.

In generalised type systems, this bijective equivalence between types
and proofs is broadened. Predicate logic is related to dependent
types systems, and higher order logic to polymorphic types. Type constructors
are identified with the usual operators and quantifiers: $\Pi$-constructor
equates to generalisation and implication, $\Sigma$-constructor to
existential quantification; tautology equates to type 1, contradiction
to void type. 

\begin{center}
\begin{tabular}{|c|c|}
\hline 
\noun{Type theory} & \noun{Logic}\tabularnewline
\hline 
\hline 
{\small{Types and terms, $a:A$}} & {\small{$a$ is a proof of $A$}}\tabularnewline
\hline 
{\small{Unit type and void type }} & {\small{Tautology and contradiction}}\tabularnewline
\hline 
{\small{Product type $A\times B$}} & {\small{ $A\land B$}}\tabularnewline
\hline 
{\small{Function type $A\rightarrow B$}} & {\small{$A\rightarrow B$}}\tabularnewline
\hline 
$\Sigma$ and $\Pi$ dependent types & Quantifiers, $\forall$ and $\exists$\tabularnewline
\hline 
Application, $(\lambda x^{A}.t:B)u^{A}\vdash(tu):B$ & \emph{Modus ponens}, $A\rightarrow B,A\vdash B$\tabularnewline
\hline 
\end{tabular}
\par\end{center}

\begin{center}
\noun{Table 3.}
\par\end{center}

\medskip{}

In sum, there is a correspondence between type systems and logical
calculus systems, and between the elements and rules of these systems,
and from an intuitionistic interpretation this correspondence serves
as a link between proof theory and type theory. This isomorphism can
be extended to the Curry\textendash{}Howard\textendash{}Lambek correspondence
if we include the isomorphism between type theories and category theory
that we have seen before, and the equivalence between CCC and intuitionistic
propositional calculus observed by Lambek (1972).

\subsection{Girard's paradox }

However, this correspondence cannot be fully maintained in certain
type systems, since Girard's paradox (Cfr. (Girard, 1972) and (Coquand,
1986)) states that a type theory cannot quantify over all propositions
and identify types and propositions at the same time. Therefore, one
of these two points has to be left aside in order to maintain the
validity of the other one. This issue will serve as a major difference
between predicative ITT and impredicative CoC. The first one identifies
types and propositions and thus leaves aside universal quantification
over propositions, whereas the second one takes away the bijective
identification between types and propositions: 
\begin{itemize}
\item In an earlier, impredicative version of ITT, types $Prop$ and $Type$
are each identified with one another. The problem arises because the
proposition $Type:Type$ is not normalizing, thus it is not a well-typed
term and it leads to a contradiction in the field of types analogous
to Burali-Forti's paradox (Cfr. (Reinhold, 1989)). As has been mentioned
before, later versions of this theory avoid this paradox by proposing
a hierarchy of predicative universes in which the bijective relationship
between propositions as types is not problematic. According to Curry-Howard
correspondence, universes can be seen both as a constructive hierarchy
of types or as a hierarchy of predicates (Palmgren, 1998). 
\item CoC prefers to maintain quantification over propositions, thus not
identifying isomorphically propositions and types. In this calculus,
the identification of each proposition with the type of its proof
is maintained, but it does not allow to identify every type with a
proposition, because in CoC there are non-propositional types. In
this way, CoC can be understood as a variation of Curry-Howard isomorphism,
since strictly speaking it does not present a real isomorphic relation
between types and predicates (Coquand, 1986), but a weaker one.
\end{itemize}
Therefore, there cannot be a single unified, normalising type theory
with the aforementioned properties, universal quantification and identification
of types and propositions.

\section{Conclusions}

To sum up, the following remarks can be made: 
\begin{itemize}
\item First of all, it is worthwhile mentioning the significance of classical
paradoxes and their roles in the foundations of mathematics. Russell's
paradox and Church's type theories are still the object of fruitful
studies. 
\item Nowadays there are two basic aspects that make the study of pure type
systems appealing. The first one is related to mathematics and the
second one, to theory of computation: Firstly, in recent decades type
theories are being studied in connection to new areas of mathematics
and more specifically to the foundations of mathematics, such as category
theory or homotopy theory. In a wider sense, it is interesting to
observe how some of the problems of generalised type theory are similar
to classical problems of Frege's logic and naïve set theory. Secondly,
type theories as computing languages are the object of intensive research
these days. Proof checkers, theorem provers and type checkers, tools
within the field of automated reasoning, are founded on the principles
of type theories. Turing incompleteness, that can sometimes be seen
as a flaw, is revealed in other contexts as an advantage, since a
Turing-incomplete type system is decidable. In this way, research
on Coq or Agda relies on a heavy study on type theories, since they
are decidably verifiable.
\end{itemize}
A few closing remarks can be made in a more general style. It can
be useful to see how, instead of a single logical theory, a plurality
of logical systems have emerged, giving new stimulus to the renewal
of logical studies. Type theories were initially created as a response
to the paradoxes of naïve set theory. This connection has never been
lost, since several contemporary authors still study these theories
as an alternative basis for the foundations of mathematics. Although
the desire of a \emph{lingua universalis}, a logical language able
to verily code the structure of reality, still persists (Cfr. (Gränstrom,
2011)), the current direction points to a different goal. Instead
of a global logical system, we can conceive several frames with invariances
between them in which several logical theories can be developed (Cfr.
Bell (1986)). The question is not which logic is the real one, but
which one we desire to use for the purposes of our research, because
each type theory and more broadly each logical system has its own
advantages and weaknesses and there is no system that can comprise
all desired logical properties into a global, unified theory. 

\bibliographystyle{plainnat}
\nocite{*}
\bibliography{AxiomathesRevised}

\begin{thebibliography}{31}
\providecommand{\natexlab}[1]{#1}
\providecommand{\url}[1]{\texttt{#1}}
\expandafter\ifx\csname urlstyle\endcsname\relax
  \providecommand{\doi}[1]{doi: #1}\else
  \providecommand{\doi}{doi: \begingroup \urlstyle{rm}\Url}\fi

\bibitem[Asperti and Longo(1991)]{asperti1991categories}
Andrea Asperti and Giuseppe Longo.
\newblock \emph{Categories, types, and structures: an introduction to category
  theory for the working computer scientist}.
\newblock MIT press, 1991.

\bibitem[Barendregt(1991)]{New1}
H.P. Barendregt.
\newblock {Introduction to generalized type systems}.
\newblock \emph{Journal of Functional Programming}, \penalty0 (1(2)):\penalty0
  125--154, 1991.

\bibitem[Barendregt(1992)]{New4}
H.P. Barendregt.
\newblock {Lambda calculi with types}.
\newblock In D.~Gabbay T.~Maibaum {S. Abramsky}, editor, \emph{{Handbook of
  Logic in Computer Science}}. Oxford Science Publications, 1992.

\bibitem[Bell(1986)]{bell1986absolute}
John~L Bell.
\newblock From absolute to local mathematics.
\newblock \emph{Synthese}, 69\penalty0 (3):\penalty0 409--426, 1986.

\bibitem[Church(1940)]{New14}
A.~Church.
\newblock {A formulation of the simple theory of types}.
\newblock \emph{Journal of Symbolic Logic}, \penalty0 (5):\penalty0 56--68,
  1940.

\bibitem[Church(1941)]{church1985calculi}
Alonzo Church.
\newblock The calculi of lambda-conversion.
\newblock \emph{Series: Annals of Mathematics Studies}, \penalty0 (6), 1941.

\bibitem[Coquand(1986)]{New3}
T.~Coquand.
\newblock {An analysis of Girard's paradox}.
\newblock \emph{Proceedings of the IEEE Symposium on Logic in Computer
  Science}, pages 227--236, 1986.

\bibitem[Curry(1934)]{curry1934functionality}
Haskell~B Curry.
\newblock Functionality in combinatory logic.
\newblock \emph{Proceedings of the National Academy of Sciences of the United
  States of America}, 20\penalty0 (11):\penalty0 584, 1934.

\bibitem[Girard(1972)]{New11}
J.Y. Girard.
\newblock \emph{{Interpretation fonctionelle et eleimination des coupures dans
  l'arithmetique d'ordre superieure}}.
\newblock PhD thesis, Universit{\'e} Paris 7, 1972.

\bibitem[Granstr{\"o}m(2011)]{New7}
J.G. Granstr{\"o}m.
\newblock \emph{{Treatise on Intuitionistic Type Theory}}.
\newblock Springer, 2011.

\bibitem[Howard(1969)]{howard1995formulae}
William~A Howard.
\newblock The formulae-as-types notion of construction. unpublished manuscript.
\newblock 1969.

\bibitem[Jacobs(1999)]{jacobs1999categorical}
Bart Jacobs.
\newblock \emph{Categorical logic and type theory}, volume 141.
\newblock Elsevier, 1999.

\bibitem[Johnstone(2003)]{johnstone2003sketches}
Peter~T Johnstone.
\newblock Sketches of an elephant: A topos theory compendium-2 volume set.
\newblock \emph{Sketches of an Elephant: A Topos Theory Compendium-2 Volume
  Set}, 1, 2003.

\bibitem[Lambek(1972)]{lambek1972deductive}
Joachim Lambek.
\newblock Deductive systems and categories iii. cartesian closed categories,
  intuitionist propositional calculus, and combinatory logic.
\newblock In \emph{Toposes, algebraic geometry and logic}, pages 57--82.
  Springer, 1972.

\bibitem[Lambek and Scott(1988)]{lambek1988introduction}
Joachim Lambek and Philip~J Scott.
\newblock \emph{Introduction to higher-order categorical logic}, volume~7.
\newblock Cambridge University Press, 1988.

\bibitem[Landin(1965)]{New2}
P.~Landin.
\newblock {Correspondence between ALGOL 60 and Church's Lambda Notation: Part
  I}.
\newblock \emph{Communications of the ACM}, 8\penalty0 (2), 1965.

\bibitem[Lawvere(1963)]{lawvere1963functorial}
F~William Lawvere.
\newblock Functorial semantics of algebraic theories.
\newblock \emph{Proceedings of the National Academy of Sciences of the United
  States of America}, 50\penalty0 (5):\penalty0 869, 1963.

\bibitem[Martin-L{\"o}f(1971)]{New15}
P.~Martin-L{\"o}f.
\newblock {A theory of types. Preprint}, 1971.

\bibitem[Martin-L{\"o}f(1975)]{New6}
P.~Martin-L{\"o}f.
\newblock {An intuitionistic theory of types: predicative part}.
\newblock \emph{Logic Colloquium '73}, pages 73--118, 1975.

\bibitem[Palmgren(1998)]{New12}
Erik Palmgren.
\newblock {On universes in type theory}.
\newblock New York, 1998.

\bibitem[Peruzzi(1991)]{peruzzi1991categories}
A~Peruzzi.
\newblock Categories and logic.
\newblock \emph{Problemi fondazionali nella teoria del significato}, pages
  137--211, 1991.

\bibitem[Ramsey(1931)]{ramsey1931foundations}
Frank~Plumpton Ramsey.
\newblock Foundations: Essays in philosophy, logic, mathematics, and economics.
\newblock 1931.

\bibitem[Ramsey(2013)]{ramsey2013foundations}
Frank~Plumpton Ramsey.
\newblock \emph{Foundations of mathematics and other logical essays}.
\newblock Routledge, 2013.

\bibitem[Reinhold(1989)]{New10}
Mark Reinhold.
\newblock {Typechecking is Undecidable When {\lq}Type{\rq} is a Type}.
\newblock Technical report, Massachusets Institute of Technology, 1989.

\bibitem[Reynolds(1974)]{reynolds1974towards}
John~C Reynolds.
\newblock Towards a theory of type structure.
\newblock In \emph{Programming Symposium}, pages 408--425. Springer, 1974.

\bibitem[Roorda(2000)]{New13}
Jan-Willem Roorda.
\newblock {Pure Type Systems for Functional Programming}.
\newblock Master's thesis, University of Utrecht, 2000.

\bibitem[Russell(1903)]{New9}
B.~Russell.
\newblock \emph{{The Principles of Mathematics}}.
\newblock Cambridge University Press, 1903.

\bibitem[Russell(1908)]{russell1908mathematical}
Bertrand Russell.
\newblock Mathematical logic as based on the theory of types.
\newblock \emph{American journal of mathematics}, 30\penalty0 (3):\penalty0
  222--262, 1908.

\bibitem[Russell(1980)]{russell1980correspondence}
Bertrand Russell.
\newblock Correspondence with frege.
\newblock \emph{Philosophical and Mathematical Correspondence}, pages 130--170,
  1980.

\bibitem[Seely(1984)]{seely1984locally}
Robert~AG Seely.
\newblock Locally cartesian closed categories and type theory.
\newblock In \emph{Mathematical proceedings of the Cambridge philosophical
  society}, volume~95, pages 33--48. Cambridge Univ Press, 1984.

\bibitem[Voevodsky(2013)]{New5}
VA~Voevodsky.
\newblock \emph{{Homotopy Type Theory}}.
\newblock Univalent Foundations Program, 2013.

\end{thebibliography}

\end{document}